\documentclass[a4paper,12pt]{article}
\usepackage{amssymb,amsmath,amsthm}
\numberwithin{equation}{section}
\usepackage{natbib}
\usepackage{bm}
\usepackage{latexsym}
\usepackage{setspace}
\usepackage{mathtools}
\usepackage{relsize}
\usepackage[margin=1in]{geometry}
\usepackage[compact]{titlesec}
\usepackage[titletoc]{appendix}
\titleformat*{\section}{\large\bfseries}
\titleformat*{\subsection}{\normalsize\bfseries}
\titleformat*{\subsubsection}{\normalsize\bfseries}
\titlespacing{\section}{0pt}{0.5cm}{0.5cm}
\titlespacing{\subsection}{0pt}{0.5cm}{0.5cm}
\usepackage[bookmarks,bookmarksnumbered,bookmarksopen,breaklinks,linktocpage]{hyperref}
\newcommand{\R}{{\mathbb R}}

\newcommand{\sgn}{\operatorname{sgn}}
\newtheorem{Theorem}{Theorem}[section]

\newtheorem{Lemma}{Lemma}[section]

\begin{document}   %general article-type format
\thispagestyle{empty}

\author{{ Jitraj Saha\footnote{\textit{ Email address:} jitraj@iitkgp.ac.in}} $ $ and $ $ Jitendra Kumar \vspace{.2cm}\\%\footnote{\textit{ Email address:} jkumar@maths.iitkgp.ernet.in}\vspace{.2cm}\\
\small \it Department of Mathematics, Indian Institute of Technology Kharagpur\\
\small \it Kharagpur-721302, India.}

%\title{Existence and uniqueness of solutions with singular coagulation and unbounded multi-fragmentation kernels}
\title{Existence and uniqueness of solutions for coagulation-fragmentation problems with singularity}

\maketitle
%\begin{spacing}{1.4}
\begin{abstract}
In this paper, existence and uniqueness of solutions to a non-linear, initial value problem is studied. In particular, we consider a special type of problem which physically represents the time evolution of particle number density resulted due to the coagulation and fragmentation process. The coagulation kernel is chosen from a huge class of functions, both singular and non-singular in nature. On the other hand, fragmentation kernel includes practically relevant non-singular unbounded functions. Moreover, both the kernels satisfy a linear growth rate of particles at infinity. The existence theorem includes lesser restrictions over the kernels as compared to the previous studies. Furthermore, strong convergence results on the sequence of functions are used to establish the existence theory.
\end{abstract}

 \vspace{0.2cm}\\
\textbf{2010 subject classification:} $45$J$05$, $34$A$34$, $45$L$10$

%%%%%%%%%%%%%%%%%%%%%%%%%%%%%%%%%%%%%% Einleitung %%%%%%%%%%%%%%%%%%%%%%%%%%%%%%%%%%%%%%%%%%%%%%%%%%%%%%%

\section{Introduction}
\label{s1}

We consider the general non-linear, initial value problem representing particle coagulation and multi-fragmentation
\begin{align}
\frac {\partial f(x,t)}{\partial t} = \frac{1}{2} \int_0^x  & K(x-y,y)f(x-y,t)f(y,t)\, dy - f(x,t) \int_0^\infty K(x,y)f(y,t)\, dy \notag \\
& + \int_x^\infty b(x,y)S(y)f(y,t)\, dy - S(x)f(x,t) \label{1.1}
\end{align}
supported with the initial data,
\begin{align}
f(x,0)=f_{0}(x)\geq 0,\quad \mbox{for all} \quad x>0. \label{1.2}
\end{align}
In equation (\ref{1.1}), the function $f(x,t)$ denotes the number density of the particles of volume $x>0$ at time $t\ge 0$. The coagulation kernel $K(x,y)$ represents the rate of coagulation (or aggregation) of particles with volume $x$ and $y$ to form larger agglomerates. The selection rate of particles of volume $y$ for fragmenting into smaller particles, is denoted by the function $S(y)$ and the distribution of daughter particles $x$ formed due to the fragmentation of $y$ is represented by the breakage function $b(x,y)$. Further, the first and the third terms of equation (\ref{1.1}) denote the formation or \emph{birth} of particles with volume $x$ in the system. On the other hand, the second and the forth terms denote the removal or \emph{death} of particles having volume $x$ from the system. In general, each of $f(x,t)$, $K(x,y)$, $S(x)$ and $b(x,y)$ are assumed to be non-negative functions. Moreover, the coagulation kernel $K(x,y)$ is considered to be a symmetric function of $x$, $y$ and the breakage function satisfies the relations
\begin{align}
\int_0^y b(x,y)\, dx \le \nu, \quad \mbox{and}\quad \int_0^y xb(x,y)\, dx = y. \label{1.3}
\end{align}
The first integral in (\ref{1.3}) represents the total number of fragments produced during the breakage event of particle volume $y$ and the second integral represents the volume conservative property of system during particle fragmentation. In this work, we consider the number of fragments to be bounded by a positive and finite quantity $\nu$.\vspace{0.2cm}\\
The study of coagulation-fragmentation (CF) equations with different forms of kernels is of great interest in the fields of engineering sciences [e.g.\ \cite{qamar2007a,qamar2007,maas2010tailoring} (chemical engineering), \cite{smikalla2011impact} (pharmaceutical), \cite{kind1999product} (food processing) etc.] as well as in mathematical sciences [\cite{stewart1991density,laurenccot2002discrete,banasiak2012analytic,ganesan2012operator}]. From theoretical point view, study on existence and uniqueness of solutions for the CF-equations is gathering high importance. In the literature, several articles by \cite{stewart1989global,stewart1990equation,stewart1990uniqueness,dacosta1995existence,dubovskii1996existence,laurenccot2000class,
Escobedo,lamb2004existence} can be found where the existence, uniqueness theory of the solutions has been established. In most of the above mentioned articles, the authors have considered a non-singular, unbounded coagulation kernel. Beside these articles, some notable works by \cite{smoluchowski1916a,smoluchowski1917,kapur1972kinetics,sastry1975similarity,Ding} can be found where the modeling of the CF-equations incorporating several physical phenomenon are done by considering singular coagulation kernels [For details, please see Appendix-\ref{Apdx1}]. In this regard, it is needed to mention that in industrial engineering sectors like pharmaceutical, food processing etc.\ dust formation of particles is unwanted due to their handling difficulty. So, the coagulation kernels are set with functions having singularity at the coordinate axes. Therefore, the dust particles tend to coagulate at a high rate to form larger agglomerates. However, the mathematical existence of solutions for the CF-equations with singular coagulation kernels has not been studied so extensively. Very recently, \cite{saha2014singular,camejo2014regular,camejo2014singular} have studied the existence theory for the problems with singular coagulation kernels. Let us briefly look through the works done in these articles. \cite{camejo2014regular} proved the existence theory using weak $L^1$ compactness theory, for a pure aggregation problem with a very restrictive class of singular coagulation kernel satisfying $$K(x,y)\le k\frac{(1+x+y)^\lambda}{(xy)^\sigma}, \quad 0\le \sigma\le 1/2,\quad 0\le \lambda-\sigma< 1.$$ On the other hand, \cite{saha2014singular} proved the existence and uniqueness of mass conserving, strong solutions of the CF-equations with almost similar type of singular coagulation kernel as used by \cite{camejo2014regular} and an unbounded class of non-singular fragmentation kernels. Later, \cite{camejo2014singular} extended the existence theory of weak solutions of the CF-equations for a wide range of singularity over the coagulation kernel. In their article, \cite{camejo2014singular} have considered $$K(x,y)\le k\frac{(1+x)^\lambda(1+y)^\lambda}{(xy)^\sigma}, \quad \sigma\ge 0,\quad 0\le \lambda-\sigma< 1$$ and the multiple fragmentation kernel is chosen to be at most linear, $$0\le S(x)\le x^\theta \quad \mbox{for}\quad \theta\in[0,1[$$ along with the following restrictions over the breakage function
$$\int_0^y b(x,y)x^{-2\sigma}\, dx\le Cy^{-2\sigma},~ \int_0^y b^q(x,y)\,dx\le B_1y^{q\tau_1},~ \int_0^y x^{-q\sigma}b^q(x,y)\, dx\le b_2y^{q\tau_2}$$ where $C,B_1,B_2$ are positive constants and $q>1$, $\tau_1,\tau_2\in [-2\sigma-\theta,1-\theta]$.%\vspace{0.2cm}\\
\subsection{State of the art}
\label{s1.1}
For a fixed $T(>0)$, we consider the strip $$\mathcal{S} := \left\{(x,t) : 0<x<\infty, 0\le t\le T \right\},$$ and define $\Upsilon_{r,\sigma} (T)$ to be the space of all continuous functions $f$ with the bounded norm
\begin{align}
\left\| f \right\|_\Upsilon := \sup_{0\le t\le T} \int_0^\infty \left( x^{r} + \frac{1}{x^{2\sigma}}\right) \left|f(x,t)\right|\, dx, \quad r\ge 1,~ \sigma\ge 0. \label{2.1}
\end{align}
It can be easily verified that $\Upsilon_{r,\sigma} (T)$ is a Banach space. Furthermore, let $\Upsilon_{r,\sigma}^+(T)$ denotes the space of all non-negative functions from $\Upsilon_{r,\sigma} (T)$. In this article, we aim to prove the existence of solutions for the CF-equation (\ref{1.1}), (\ref{1.2}) under the following assumptions over the coagulation and fragmentation kernels;
\begin{enumerate}
\item[(A1)] $K(x,y)\le \displaystyle{k\frac{(1+x+y)^\lambda}{(xy)^\sigma}}$, \quad where \quad $\sigma\ge 0$,\quad $0\le \lambda-\sigma\le 1$,
\item[(A2)] $S(x)\le s_0 x^\alpha$, where $0\le \alpha \le \displaystyle{\min\left[r,1+2\sigma \right]}$, \quad $r\ge 1$, and
\item[(A3)] $\displaystyle{\int_0^y x^{-2\sigma}b(x,y)\, dx} \le \bar{N}y^{-2\sigma}$ and for any $y$ satisfying $y\ge x$, $\displaystyle{\sup_{x\in \left[ x_1,x_2\right]}} xb(x,y)\le \bar{b}$, for all $0<x_1<x_2<\infty$.
\end{enumerate}
Here, $k$, $s_0$, $\bar{N}$ and $\bar{b}$ are considered to be positive constants.\vspace{0.2cm}\\
Substituting $\sigma=0$ and $\lambda=1$ in $(A1)$, we obtain the coagulation kernel considered by \cite{dubovskii1996existence}. Similarly, $K(x,y)$ of \cite{camejo2014regular} can be obtained from $(A1)$, whenever $0\le \sigma\le 1/2$. In this regard, using the inequality $(x+y)^p \le  2^{p}\left(x^p+y^p\right)$ for all $x,y,p\ge 0$,
%$$\left(x+y\right)^p \le \left(x^p+y^p\right), \quad \mbox{if } 0\le p\le 1, \quad \mbox{and} \quad (x+y)^p \le  2^{p-1}\left(x^p+y^p\right), \quad \mbox{if } p \ge 1,$$
%$$
%\begin{array}{lll}
%\left(x+y\right)^p & \le \left(x^p+y^p\right), & \quad \mbox{if } 0\le p\le 1,\\
%(x+y)^p & \le  2^{p-1}\left(x^p+y^p\right), & \quad \mbox{if } p \ge 1.
%\end{array}
%$$
the coagulation kernel of \cite{camejo2014singular} can be written as
\begin{align*}
K(x,y) \le k\frac{(1+x+y+xy)^\lambda}{(xy)^\sigma} \le 2^\lambda k\left[\frac{(1+x+y)^\lambda}{(xy)^\sigma} +(xy)^{\lambda-\sigma}\right].
\end{align*}
Thus, the above relation indicates that $K(x,y)$ of \cite{camejo2014singular} includes an additional non-singular class of functions in the form $(xy)^{\lambda-\sigma}$ with $0\le \lambda-\sigma\le 1$, than the class considered in $(A1)$. However, coagulation kernel involving the above class of non-singular, unbounded functions has already been dealt in our previous article \cite{saha2014singular}. So in the present study, we emphasis on the strong singularity conditions over $x$ and $y$ by widening the working range of $\sigma$. Thus, the coagulation kernel of \cite{camejo2014singular} and $(A1)$ can be treated to be equivalent when the class of singular functions are taken into account. Moreover, with the above mentioned bounds over the coagulation term, we are able to include the generalized form of Brownian motion kernel [\cite{smoluchowski1917}] and granulation kernel [\cite{kapur1972kinetics}] into our consideration.\vspace{0.2cm}\\
Besides the coagulation kernel, the selection function considered in this work $(A2)$, includes a wider class of unbounded function compared to that of \cite{camejo2014singular}. Furthermore, we impose lesser restrictions over the breakage functions $(A3)$ than that of \cite{camejo2014singular}. In this regard, let us take an example of the following breakage function $$b(x,y) = (\tau+2)\frac{x^\tau}{y^{\tau+1}}, \quad \mbox{with} \quad \tau+1>0,$$ which are well known in the literature as the `power law kernels' [\cite{banasiak2013strong}]. Let $\sigma$ is so chosen that $\left(\tau+1\right)>2\sigma$, then
\begin{align*}
\int_0^y x^{-2\sigma}b(x,y)\, dx = \frac{(\tau+2)}{y^{\tau+1}}\int_0^y x^{\tau-2\sigma} \, dx = \frac{(\tau+2)}{(\tau-2\sigma+1)}\frac{ y^{\tau-2\sigma+1}}{y^{\tau+1}} = \frac{(\tau+2)}{(\tau-2\sigma+1)}y^{-2\sigma}.
\end{align*}
Also, for $0<x_1\le x\le x_2<\infty$, where $x_1$, $x_2$ are finite and $x\le y$, we get $$\sup_{x\in[x_1,x_2]} xb(x,y) \le (\tau+2)\left(\frac{x}{y}\right)^{\tau+1}\le (\tau+2).$$
Thus, the condition $(A3)$ is a realistic assumption on $b(x,y)$ which is satisfied by the power law kernels.\vspace{0.2cm}\\
Our existence result is motivated upon the work of \cite{dubovskii1996existence}. In the next section, we state and prove the existence theorem. Since, strong convergence results are used to prove the theorem so, the existence of strong solutions to the initial value problem (IVP) (\ref{1.1}), (\ref{1.2}) is obtained. In section \ref{s3}, the uniqueness property of the solutions is studied.

\section{Existence theorem}
\label{s2}

%%Let for a fixed $T(>0)$, we define the space $$\mathcal{S} := \left\{(x,t) : 0<x<\infty, 0\le t\le T \right\}.$$ Further, consider that for $ r\ge 1$, $\sigma\ge 0$ and $\Upsilon_{r,\sigma} (T)$ denotes the space of all continuous functions $f(x,t)$ from $\mathcal{C}\left( ]0,\infty[,[0,T]\right)$, where $0<x<\infty$ and $0\le t\le T$. Therefore, $\Upsilon_{r,\sigma} (T)$ becomes a Banach space when equipped with the norm
%%\begin{align}
%%\left\| f \right\|_\Upsilon := \sup_{0\le t\le T} \int_0^\infty \left( x^{r} + \frac{1}{x^{2\sigma}}\right) \left|f(x,t)\right|\, dx < \infty. \label{2.1}
%%\end{align}
%%Let the cone of all non-negative functions in $\Upsilon_{r,\sigma}(T)$ be denoted by $\Upsilon_{r,\sigma}^+(T)$.

\begin{Theorem}
\label{th1}
Let the functions $K(x,y)$, $b(x,y)$ and $S(x)$ be non-negative and continuous for all $x,y\in]0,\infty[$ and satisfy the condition $(A1)$, $(A2)$ and $(A3)$. If the initial data $f_0(x)$ belongs to $\Upsilon_{r,\sigma}^+(0)$, then the IVP (\ref{1.1}), (\ref{1.2}) has at least one solution in $\Upsilon_{r,\sigma}^+(T)$.
\end{Theorem}

%%The proof has two cases.
%%\subsection*{Case$-1$:} Let us first consider that there exists a fixed $R(>0)$, such that for each $t\in[0,T]$ the intervals $\left[\frac{1}{R},R\right]\times\left[\frac{1}{R},R\right]$ and $\left[0,R\right]$, are the compact support of the kernels $K(x,y)$ and $S(x)$, respectively. Then, for $x,y>R$ and $x,y<\frac{1}{R}$ the solution of (\ref{1.1}), (\ref{1.2}) is given by
%%\begin{align}
%%f(x,t) = f_0(x). \label{2.2}
%%\end{align}
%%The relation (\ref{2.2}) gives an estimate of the solution function outside the compact domain. Thus we get that for problems with compactly supported kernels, the \emph{tails} of the solution $f(x,t)$ actually does not change with time and coincides with the \emph{tails} of the initial data $f_0(x)$. The existence and uniqueness of the continuous solutions within the support can easily be obtained using Banach fixed point theorem over a suitably defined invariant closed ball. The proof of this part bears similarity to the work of \cite{dubovskii1994mathematical} and thus can easily be followed from there.
%
%%\subsection*{Case$-2$:}
%%Let us consider that the kernels $K(x,y)$ and $S(x)$ do not possess any compact support.
\begin{proof}
Recalling the assumptions $(A1)$ and $(A2)$, we observe that the kernels $K(x,y)$ and $S(x)$ exhibit unbounded rates. Therefore, we first truncate both the coagulation and selection function in order to construct the following sequence of continuous kernels $\left\{ K_n,S_n \right\}_{n=1}^\infty$ with compact support for each $n\ge1$, from the class defined in $(A1)$ and $(A2)$

$$K_{n}(x,y)\left\{\begin{array}{ll}
     =K(x,y), &\mbox{when } x,y\ge \frac{1}{n}, \mbox{ and } x+y \le n, ,\\
     \le K(x,y), &\mbox{elsewhere,}
     \end{array}
\right.$$
and
$$S_{n}(x)\left\{\begin{array}{ll}
     =S(x), &\mbox{when $0\le x\le n$,}\\
     \le S(x), &\mbox{elsewhere.}
     \end{array}
\right.$$

Therefore, for the above defined `cut-off' kernels $K_n$ and $S_n$ the IVP (\ref{1.1}) is written as
\begin{align}
\frac{\partial f_n(x,t)}{\partial t} = & \frac{1}{2}\int_0^x K_n(x-y,y)f_n(x-y,t)f_n(y,t)\, dy - f_n(x,t)\int_0^\infty K_n(x,y)f_n(y,t)\, dy \notag \\
& \quad + \int_x^\infty b(x,y)S_n(y)f_n(y,t)\, dy - S_n(x)f_n(x,t), \label{2.4}
\end{align}
with the truncated initial data
\begin{align}
f_{n}(x,0)\left\{\begin{array}{ll}
     =f_0(x), &\mbox{when $0\le x\le n$,}\\
     \le f_0(x), &\mbox{elsewhere.}
     \end{array}
\right.
\end{align}

Now for each $n\ge 1$, the existence of continuous non-negative solutions $f_n$ for the problem (\ref{2.4}) can be established by following the works of \cite{stewart1989global,camejo2014singular}. Thus a sequence of continuous, non-negative solutions $\{f_n\}_{n=1}^\infty$ to the problem (\ref{1.1}) is obtained which satisfy
\begin{align}
f_n (x,t)\in\Upsilon_{r,\sigma}^+ (T), \quad \mbox{for each}\quad n\ge 1. \label{2.3}
\end{align}
%Referring to Case$-1$ we get that for each $n\ge 1$, the truncated kernels $K_n$ and $S_n $ generates a sequence of continuous, non-negative solutions $\{f_n\}_{n=1}^\infty$ to the problem (\ref{1.1}). Moreover,

Now corresponding to the solutions $f_n$, we define the $(r+2)-$number of moments as
\begin{align}
M_{i,n}(t) = \int_0^\infty x^i f_n(x,t)\, dx,\quad  i=-2\sigma, 0,1,2,\ldots,r, \quad \mbox{and} \quad n\ge 1, \label{2.5}
\end{align}
and establish the boundedness of these $(r+2)$ moments.
\subsection{Boundedness of the moments $M_{i,n}(t)$:}
For $i=1$, we have
\begin{align*}
\frac{dM_{1,n}(t)}{dt} & = \int_0^\infty x\frac{\partial f_n(x,t)}{\partial t}\, dx.
\end{align*}
Since $K_n$ and $S_n$ have compact support, all the integrals appeared after computing right hand side of the above equation are finite and they cancel out. Hence, for all $n\ge 1$ and $0\le t\le T$ we get a constant independent of $n$, such that
\begin{align}
M_{1,n}(t) = \int_0^\infty xf_0(x)\, dx = \bar{M}_1. \label{2.6}
\end{align}
Now proceed to check the behavior of the moment $M_{-2\sigma,n}$. Therefore, multiplying with the weight $x^{-2\sigma}$ and integrating, equation (\ref{2.4}) is written as
\begin{align*}
\frac{dM_{-2\sigma,n} (t)}{dt} = & \frac{1}{2}\int_0^\infty x^{-2\sigma}\int_0^x K_n(x-y,y)f_n(x-y,t)f_n(y,t)\, dy\, dx \notag \\
& - \int_0^\infty x^{-2\sigma}f_n(x,t)\int_0^\infty K_n(x,y)f_n(y,t)\, dy\,dx \notag \\
& + \int_0^\infty x^{-2\sigma}\left[\int_x^\infty b(x,y)S_n(y)f_n(y,t)\, dy - S_n(x)f_n(x,t)\right]\, dx. %\label{2.7}
\end{align*}
Let us now perform some mathematical computations,
\begin{enumerate}
\item[(i)] change the order of integration of the first integral and third integral in the r.h.s.\ of the above equation.
\item[(ii)] substitute $x-y = \bar{x}$, $y = \bar{y}$ and again replace $\bar{x} = x$ and $\bar{y} = y$, to get
\end{enumerate}
\begin{align}
\frac{dM_{-2\sigma,n} (t)}{dt} = & \frac{1}{2}\int_0^\infty \int_0^\infty \left[(x+y)^{-2\sigma} - x^{-2\sigma} -y^{-2\sigma}\right] K_n(x,y)f_n(y,t)\, dx\, dy \notag \\
& + \int_0^\infty \int_0^y x^{-2\sigma} b(x,y)S_n(y)f_n(y,t)\, dx\, dy - \int_0^\infty x^{-2\sigma}S_n(x)f_n(x,t)\, dx.\label{2.7}
\end{align}
For $\sigma\ge 0$ and $x,y>0$, we have the inequality
\begin{align*}
(x+y)^{-2\sigma}\le x^{-2\sigma} + y^{-2\sigma}. %\label{2.9}
\end{align*}
Therefore, using the positivity of the functions $S_n$, $f_n$ and the above inequality in the relation \eqref{2.7}, we get
\begin{align*}
\frac{dM_{-2\sigma,n} (t)}{dt} \le \int_0^\infty \int_0^y x^{-2\sigma} b(x,y)S_n(y)f_n(y,t)\, dx\, dy - \int_0^\infty x^{-2\sigma}S_n(x)f_n(x,t)\, dx.
\end{align*}
Using  $(A3)$ in the above relation, we obtain
\begin{align}
\frac{dM_{-2\sigma,n} (t)}{dt} \le & \left(\bar{N}-1\right) \int_0^\infty y^{-2\sigma} S_n(y)f_n(y,t)\, dy \notag \\
\le & s_0\left(\bar{N}-1\right) \int_0^\infty y^{\alpha-2\sigma} f_n(y,t)\, dy \notag \\
= & s_0\left(\bar{N}-1\right) \left[ \int_0^1 y^{\alpha-2\sigma} f_n(y,t)\, dy + \int_1^\infty y^{\alpha-2\sigma} f_n(y,t)\, dy\right].\label{2.8}
\end{align}
From the condition $(A2)$, we have $(\alpha-2\sigma)\le 1$.
\begin{enumerate}
\item[(i)] If $0\le (\alpha-2\sigma)\le 1$, then $y\le 1$ implies $y^{\alpha-2\sigma}\le 1 \le y^{-2\sigma}$, and $y> 1$ implies $y^{\alpha-2\sigma}< y$.
\item[(ii)] If $(\alpha-2\sigma)<0$, then $y\le 1$ implies $1\le y^{\alpha-2\sigma}\le y^{-2\sigma}$, and $y> 1$ implies $y^{\alpha-2\sigma}< 1 < y$.
\end{enumerate}
Therefore, the relation (\ref{2.8}) is written as
\begin{align*}
\frac{dM_{-2\sigma,n} (t)}{dt} \le & s_0\left(\bar{N}-1\right) \left[ \int_0^1 y^{-2\sigma} f_n(y,t)\, dy + \int_1^\infty yf_n(y,t)\, dy \right] \notag \\
\le & s_0\left(\bar{N}-1\right) \left[M_{-2\sigma,n}(t) + \bar{M}_1\right]. %\label{2.11}
\end{align*}
From the above relation, we find that for $t\in [0, T]$
\begin{align}
M_{-2\sigma,n} (t) \le \left[s_0\left(\bar{N}-1\right) +1\right]\bar{M}_1\exp\left(s_0\bar{N}T\right) =: \bar{M}_{-2\sigma}. \label{2.9}
\end{align}
Here, $\bar{M}_{-2\sigma}$ is a constant term, independent of $n$. Thus the uniform boundedness of the truncated moment $M_{-2\sigma,n}$ is obtained. We now move forward to obtain the uniform boundedness of the other truncated moments of (\ref{2.5}) for $i=2,\ldots,r$. Integrating the equation (\ref{2.4}) with the weight $x^2$ and using the bounds over the kernels and proceeding as before, we obtain
\begin{align*}
\frac{dM_{2,n}(t)}{dt} \le & k\int_0^\infty \int_0^\infty (xy)^{1-\sigma}(1+x+y)^{\lambda}f_n(x,t)f_n(y,t) \, dy \, dx \\
\le & k\int_0^\infty \int_0^\infty \left[ (xy)^{1-\sigma} + x^{1+\lambda-\sigma}y^{1-\sigma} + x^{1-\sigma}y^{1+\lambda-\sigma}\right]f_n(x,t)f_n(y,t) \, dy \, dx \\
= & k \int_0^1\int_0^1 \left[ (xy)^{1-\sigma} + x^{1+\lambda-\sigma}y^{1-\sigma} + x^{1-\sigma}y^{1+\lambda-\sigma}\right]f_n(x,t)f_n(y,t) \, dy \, dx \\
& + k \int_0^1\int_1^\infty \left[ (xy)^{1-\sigma} + x^{1+\lambda-\sigma}y^{1-\sigma} + x^{1-\sigma}y^{1+\lambda-\sigma}\right]f_n(x,t)f_n(y,t) \, dy \, dx \\
& + k \int_0^\infty\int_0^1 \left[ (xy)^{1-\sigma} + x^{1+\lambda-\sigma}y^{1-\sigma} + x^{1-\sigma}y^{1+\lambda-\sigma}\right]f_n(x,t)f_n(y,t) \, dy \, dx \\
& + k \int_1^\infty\int_1^\infty \left[ (xy)^{1-\sigma} + x^{1+\lambda-\sigma}y^{1-\sigma} + x^{1-\sigma}y^{1+\lambda-\sigma}\right]f_n(x,t)f_n(y,t) \, dy \, dx.
\end{align*}
If $\left(1-\sigma \right)>0$, then obviously $0< \left(1-\sigma \right) \le 1$ and
\begin{align*}
\frac{dM_{2,n}(t)}{dt} \le  k\left[ 3\bar{M}_1^2 + 4\bar{M}_1\bar{M}_{-2\sigma} + \bar{M}^2_{-2\sigma} \right] + 4k\bar{M}_1M_{2,n}(t).
\end{align*}
If $\left(1-\sigma \right)\le 0$, then
\begin{align*}
\frac{dM_{2,n}(t)}{dt} \le k\left[ 3\bar{M}_1^2 + 4\bar{M}_1\bar{M}_{-2\sigma} + \bar{M}^2_{-2\sigma} \right] + 2k\left[ \bar{M}_1 + \bar{M}_{-2\sigma}\right]M_{2,n}(t).
\end{align*}
Recalling equations (\ref{2.6}) and (\ref{2.9}), we have $\bar{M}_1$ and $\bar{M}_{-2\sigma}$ to be constants independent of $n$ respectively. Thus from the above inequality, we can easily obtain
\begin{align}
M_{2,n} (t)\le \bar{M}_2, \mbox{ for all } 0\le t\le T, ~n\ge 1, \label{2.10}
\end{align}
where the constant $\bar{M}_2$ is independent of $n$.\vspace{0.2cm}\\
Similarly, for $r\ge3$ and integrating (\ref{2.4}) with the weight $x^3$, we get
\begin{align*}
\frac{dM_{3,n}(t)}{dt} \le \frac{3k}{2}\int_0^\infty \int_0^\infty (xy)^{1-\sigma}(x+y)(1+x+y)^{\lambda}f_n(x,t)f_n(y,t) \, dy \, dx.
\end{align*}
Simplifying, we get
\begin{align*}
\frac{dM_{3,n}(t)}{dt} \le \left \{ \begin{array}{ll}
3k\left[3\bar{M}_1\bar{M}_{-2\sigma} + 2\bar{M}_1^2+\bar{M}_2^2+3\bar{M}_1\bar{M}_2+\bar{M}_2\bar{M}_{-2\sigma} + \left(\bar{M}_1 +\bar{M}_{-2\sigma}\right)M_{3,n}\right], & \vspace{0.2cm}\\
\quad \mbox{whenever $0\le \sigma <1,$} \vspace{0.2cm}\\
3k\left[\bar{M}_{-2\sigma}^2 +5\bar{M}_1\bar{M}_{-2\sigma}+\bar{M}_1^2+\bar{M}_2\bar{M}_{-2\sigma}+\bar{M}_1\bar{M}_2 + \bar{M}_2^2 + \left(\bar{M}_1+\bar{M}_{-2\sigma}\right)M_{3,n} \right], & \vspace{0.2cm}\\
\quad \mbox{whenever $1\le \sigma <2,$} \vspace{0.2cm}\\
3k\left[\bar{M}_{-2\sigma}^2+ 4\bar{M}_1\bar{M}_{-2\sigma} +3\bar{M}_1^2 + \bar{M}_2\bar{M}_{-2\sigma} + \bar{M}_1\bar{M}_2 + \left(\bar{M}_1+\bar{M}_{-2\sigma}\right)M_{3,n}\right], & \vspace{0.2cm}\\
\quad \mbox{whenever $\sigma \ge 2.$}
\end{array}
\right.
\end{align*}
Hence a constant $\bar{M}_3$ independent of $n$ can be obtained for all $0\le t\le T$ and $n\ge 1$, such that
\begin{align}
M_{3,n} (t)\le \bar{M}_3. \label{2.11}
\end{align}
So, for $i>3$ we can proceed further in a similar way to obtain the uniform boundedness of the truncated moments, that is, $M_{i,n}(t)\le \bar{M}_i$ where, $i=1,2,\ldots,r$, $n\ge 1$ and $0\le t\le T$. Next, to get the uniform boundedness of $M_{0,n}(t)$, we proceed as below
\begin{align*}
\frac{dM_{0,n}(t)}{dt} = \int_0^\infty & \left[ \frac{1}{2} \int_0^x K_n(x-y,y)f_n(x-y,t)f_n(y,t)\, dy - \int_0^\infty K_n(x,y)f_n(x,t)f_n(y,t)\, dy\right. \notag \\
& \quad + \left. \int_x^\infty b(x,y)S_n(y)f_n(y,t)-S_n(x)f_n(x,t)\right] \, dx. %\label{2.15}
\end{align*}
Using the positivity of the functions $b$, $S_n$, $f_n$ and performing a similar computation as done to obtain (\ref{2.7}), we get
\begin{align*}
\frac{dM_{0,n}(t)}{dt} & \le \int_0^\infty\int_0^y b(x,y)S_n(y)f_n(y,t)\, dx\, dy-\int_0^\infty S_n(x)f_n(x,t)\, dx \\
& \le \left(\nu-1\right) s_0 \int_0^\infty y^\alpha f_n(y,t)\, dy \\
& = \left(\nu-1\right) s_0 \left[ \int_0^1 y^\alpha f_n(y,t)\, dy + \int_1^\infty y^\alpha f_n(y,t)\, dy \right] \\
& \le \left(\nu-1\right) s_0 \left[ \int_0^1 f_n(y,t)\, dy + \int_1^\infty y^{r}f_n(y,t)\, dy \right] \\
& \le \left(\nu-1\right) s_0 \left[ M_{0,n}(t) + \bar{M}_{r}\right].
\end{align*}
Hence for all $0\le t\le T$ and $n\ge 1$, we have
\begin{align}
M_{0,n}(t) \le \bar{M}_0, \label{2.12}
\end{align}
where $\bar{M}_0$ is a constant independent of $n$.\vspace{0.2cm}\\
Therefore, combining all these above relations we obtain constants $\bar{M}_i$ (independent of $n$), such that
\begin{align}
M_{i,n}(t) \le \bar{M}_i, \mbox{ for all } i=-2\sigma,0,1,\ldots,r,~t\in [0,T] \mbox{ and } n\ge 1. \label{2.13}
\end{align}
We now state and prove the following lemma.
\begin{Lemma}
\label{l1}
The sequence of solutions $\left\{ f_n\right\}_{n=1}^\infty $ is relatively compact in the uniform-convergent space of continuous functions over a compact rectangular subset of $\mathcal{S}$.
\end{Lemma}

\begin{proof}
We define a compact rectangular strip of $\mathcal{S}$ as $$\bar{\mathcal{S}} := \left\{ (x,t) : \frac{1}{X} \le x\le X,~0\le t\le T \right\},$$ where $X(\gg 1)$ is any finite real number. This lemma is proved in the following steps,
\begin{enumerate}
\item uniform boundedness of the sequence $\{ f_n \}$ is obtained over $\bar{\mathcal{S}}$, and
\item then the equi-continuity of the sequence $\{ f_n \}$ with respect to both the volume variable $x$ and time variable $t$ is established over a subset of $\bar{\mathcal{S}}$.
\end{enumerate}

\subsection{Uniform boundedness of $\{f_n\}_{n=1}^\infty$:}
\label{s2.1}
In order to deal with the singularity present in the kernel $K_n$, we consider the following mathematically equivalent form of equation (\ref{2.4})
\begin{align}
x\frac{\partial f_n(x,t)}{\partial t} = -\frac{\partial \mathcal{C}(f_n)(x,t)}{\partial x} + \frac{\partial \mathcal{F}(f_n)(x,t)}{\partial x}, \label{2.14}
\end{align}
where
\begin{align*}
\mathcal{C}(f_n)(x,t) := \int_0^x \int_{x-u}^\infty uK_n(u,v)f_n(u,t)f_n(v,t)\, dv\, du,
\end{align*}
and
\begin{align*}
\mathcal{F}(f_n)(x,t) := \int_{x}^\infty \int_0^x  ub(u,v)S_n(v)f_n(v,t)\, du\, dv.
\end{align*}
Using the positivity of $K_n$, $f_n$ and the fact that $f_n\in\Upsilon_{r,\sigma}^+ (T)$, we can estimate that
\begin{align*}
x\frac{\partial f_n(x,t)}{\partial t} \le \frac{\partial \mathcal{F}(f_n)(x,t)}{\partial x}.
\end{align*}
Furthermore, recalling assumption $(A2)$ and using the fact that $x\le v$ and $u\le x$ in the term $\mathcal{F}$, we can write
\begin{align*}
x\frac{\partial f_n(x,t)}{\partial t} \le \frac{\partial \left[ x\bar{b} s_0\bar{M}_{\alpha}\right]}{\partial x} = \bar{b} s_0\bar{M}_{\alpha}.
\end{align*}
Considering $\bar{f}_0 := \displaystyle{\max_{\frac{1}{X}\le x\le X}}~ f_n(x,0)$, we have
\begin{align}
f_n (x,t) \le \bar{f}_0 + X\bar{b} s_0\bar{M}_{\alpha} =: L \quad \mbox{(say)}, \label{2.15}
\end{align}
for all $(x,t)\in \bar{\mathcal{S}}$. This $L$ is a constant independent of $n$, and thus we establish the uniform bondedness of $f_n$ over the rectangular strip $\bar{\mathcal{S}}$.

\subsection{Equicontinuity w.r.t.\ volume$-$variable $x$:}
\label{s2.2}
Let $\frac{1}{X}\le x < x'\le X$. Therefore, to establish the equicontinuity of $f_n$ with respect $x$, we need to show that, corresponding to arbitrary $\epsilon >0$, there exists an $\delta (\epsilon)$ such that
\begin{align*}
|x'-x| < \delta(\epsilon) \quad\mbox{implies}\quad \left|f_n(x',t)-f_n(x,t)\right| < \epsilon.
\end{align*}
According to the construction, each of the kernels $K_n$ and $S_n$ are continuous over the closed rectangles $\left[\frac{1}{X},X\right]\times \left[z_1,z_2\right]$ and $\left[\frac{1}{X},X\right]$, respectively. Here, we consider $0<z_1\le y\le z_2<\infty$ and these $z_1$ and $z_2$ will be determined in the following study. Furthermore, $b(x,y)$ is continuous over $[\frac{1}{X},X]\times [z_1,z_2]$. Therefore, whenever $0<\delta(\epsilon)< \epsilon$ the following relations hold true;
\begin{align*}
\sup_{|x'-x|<\delta} |f_0(x')-f_0(x)|<\epsilon, & \quad \quad \sup_{|x'-x|<\delta} \left|K_n(x',y)-K_n(x,y)\right|<\epsilon,\\
\sup_{|x'-x|<\delta} \left|S_n(x')-S_n(x)\right|<\epsilon, & \quad \quad \sup_{|x'-x|<\delta} \left|b(x',y)-b(x,y)\right|<\epsilon.
\end{align*}
Now for each $n\ge 1$, we write
\begin{align}
\left|f_n(x',t)-f_n(x,t)\right| \le & \left|f_0(x')- f_0(x)\right| \notag \\
& + \int_0^t \left[\frac{1}{2}\underbrace{\int_x^{x'} K_n(x'-y,y)f_n(x'-y,s)f_n(y,s)\, dy}_{I_1}\right. \notag\\
& + \frac{1}{2}\underbrace{\int_0^x \left| K_n(x'-y,y) - K_n(x-y,y) \right|f_n(x'-y,s)f_n(y,s)\, dy}_{I_2} \notag\\
& + \frac{1}{2}\underbrace{\int_0^x K_n(x-y,y)\left|f_n(x'-y,s) - f_n(x-y,s)\right|f_n(y,s)\, dy}_{I_3} \notag\\
& + \left|f_n(x',s)-f_n(x,s)\right|\underbrace{\int_0^\infty K_n(x',y)f_n(y,s)\, dy}_{I_4} \notag\\
& + f_n(x,s)\underbrace{\int_0^\infty \left|K_n(x',y) - K_n(x,y)\right|f_n(y,s)\, dy}_{I_5} \notag\\
& + \underbrace{\int_{x}^\infty \left|b(x',y) - b(x,y)\right|S_n(y)f_n(y,s)\, dy}_{I_6} \notag\\
& + \underbrace{\int_x^{x'} b(x',y)S_n(y)f_n(y,s)\, dy}_{I_7} + \underbrace{\left|S_n(x') - S_n(x)\right|f_n(x',s)}_{I_8} \notag\\
& + \left. \underbrace{S_n(x)\left|f_n(x',s) - f_n(x,s)\right|}_{I_9}\right] \, ds. \label{2.20}
\end{align}
Let us define
$$\zeta_n(t) :=\sup_{|x'-x|<\delta} \left|f_n(x',t)-f_n(x,t)\right|, \quad \frac{1}{X}\le x,x'\le X.$$
We first estimate the integral $I_5$ having infinite range as follows,
\begin{align*}
I_5 \le & \underbrace{\int_0^{z_1} \left|K_n(x',y)-K_n(x,y)\right|f_n(y,s)\, dy}_{J_1} + \underbrace{\int_{z_1}^{z_2} \left|K_n(x',y)-K_n(x,y)\right|f_n(y,s)\, dy}_{J_2} \\
& + \underbrace{\int_{z_2}^\infty \left|K_n(x',y)-K_n(x,y)\right|f_n(y,s)\, dy}_{J_3}.
\end{align*}
Now, the integral $J_1$ is written as
\begin{align*}
J_1 & \le 2k\int_0^{z_1} \frac{(1+x+y)^\lambda}{(xy)^\sigma}f_n(y,s)\, dy \\
& \le 2^\lambda k\int_0^{z_1} \frac{1+x^\lambda+y^\lambda}{(xy)^\sigma}f_n(y,s)\, dy \\
& \le 2^\lambda k\frac{1+x^\lambda}{x^\sigma}\int_0^{z_1} \frac{y^{\sigma}}{y^{2\sigma}} f_n(y,s)\, dy + 2^\lambda k\frac{1}{x^\sigma}\int_0^{z_1} \frac{y^{\lambda+\sigma}}{y^{2\sigma}} f_n(y,s)\, dy \\
& \le 2^\lambda k(1+X)X^\sigma z_1^{\sigma}\int_0^{z_1} \frac{1}{y^{2\sigma}} f_n(y,s)\, dy + 2^\lambda kX^\sigma z_1^{\lambda+\sigma}\int_0^{z_1} \frac{1}{y^{2\sigma}}f_n(y,s)\, dy \\
& \le 2^\lambda k(1+X)X^\sigma \left[ \bar{M}_{-2\sigma} z_1^{\sigma} + \bar{M}_{-2\sigma}z_1^{\lambda+\sigma} \right].
\end{align*}
We choose $z_1$, such that $z_1^{\sigma}\bar{M}_{-2\sigma}\le \epsilon$ and $z_1^{\lambda+\sigma}\bar{M}_{-2\sigma}\le \epsilon$. Therefore,
\begin{align*}
J_1 \le 2^{\lambda+1}k(1+X)X^\sigma\epsilon. %\label{2.23}
\end{align*}
Since, $K_n$ is equicontinuous over $\left[\frac{1}{X},X\right]\times \left[z_1,z_2\right]$, therefore we have
\begin{align*}
J_2 \le \epsilon\bar{M}_0. %\label{2.24}
\end{align*}
Next, the integral $J_3$ gives
\begin{align*}
J_3 & \le k\int_{z_2}^\infty \frac{(1+y)^\lambda}{y^\sigma}\left|\frac{(1+x')^\lambda}{(x')^\sigma}-\frac{(1+x)^\lambda}{x^\sigma}\right|f_n(y,s)\, dy \notag \\
& \le k \int_{z_2}^\infty \frac{(1+y^\lambda)}{y^\sigma}\left|\left( 1+ \frac{1}{x'} \right)^\sigma(1+x')^{\lambda -\sigma} - \left( 1+ \frac{1}{x} \right)^\sigma(1+x)^{\lambda -\sigma}\right|f_n(y,s)\, dy \notag \\
& \le k\delta\left( 1+ X\right) \int_{z_2}^\infty \left(y^{\lambda-\sigma} +\frac{1}{y^\sigma} \right)f_n(y,s)\, dy.
\end{align*}
Let $\phi (x)$ be non-negative, measurable and $\psi (x)$ is a positive and non-decreasing function for $x>0$, then we have
\begin{align}
\int_z^\infty \phi(x)\, dx \le \frac{1}{\psi (z)}\int_0^\infty \phi(x)\psi(x)\, dx, \quad z>0, \label{2.21}
\end{align}
when the integrals exists and are finite. Now in (\ref{2.21}), we first consider $\phi(x)=x^{\lambda-\sigma}f_n(x)$, $\psi(x)=x$ and then $\phi(x)=f_n(x)$, $\psi(x)=x$ respectively, to get
\begin{align*}
\int_{z_2}^\infty y^{\lambda-\sigma}f_n(y,s)\, dy \le \frac{1}{z_2}\bar{M}_2,\quad \mbox{and} \quad \int_{z_2}^\infty \frac{1}{y^\sigma}f_n(y,s)\, dy\le \frac{1}{z_2^\sigma} \int_{z_2}^\infty f_n(y,s)\, dy \le \frac{1}{z_2^{1+\sigma}}\bar{M}_1.
\end{align*}
Using the above relations, we get
\begin{align*}
J_3 \le k\delta\left( 1+ X\right)\left[ \frac{1}{z_2^{1+\sigma}}\bar{M}_1 + \frac{1}{z_2}\bar{M}_2 \right]\le 2k\delta(1+X)\epsilon, %\label{2.27}
\end{align*}
where $z_2$ is so chosen that $\frac{1}{z_2^{1+\sigma}}\bar{M}_1 < \epsilon$ and $\frac{1}{z_2}\bar{M}_2 < \epsilon$. Combining the estimates of $J_1$, $J_2$ and $J_3$, we get
\begin{align*}
I_5 \le \left[ \bar{M}_0 + 2^{\lambda+1}k(1+X)X^\sigma + 2k\delta(1+X)\right]\epsilon.%\label{1.39}
\end{align*}
The term $I_4$ is taken care in a similar manner as $I_5$ to obtain a constant $\mathcal{B}_1$, such that
\begin{align*}
\int_0^t |f_n(x',s) - f_n(x,s)|I_4 & \le \mathcal{B}_1\int_0^t \zeta_n(s)\, ds, %\label{1.40}
\end{align*}
In order to estimate $I_6$, we write
\begin{align*}
I_6 = & \int_{x'}^x \left|b(x',y)-b(x,y)\right|S_n(y)f_n(y,s)\, dy + \int_{x'}^\infty \left|b(x',y)-b(x,y)\right|S_n(y)f_n(y,s)\, dy \\
\le & \epsilon s_0\bar{M}_{\lceil \alpha \rceil} + 2\int_{x'}^\infty b(x',y)S_n(y)f_n(y,s)\, dy.
\end{align*}
Recalling the relation (\ref{2.21}), we set $\phi(y) = b(x',y)S_n(y)f_n(y)$ and $\psi(y) = y$ and get
\begin{align*}
I_6 \le & \epsilon s_0 \bar{M}_{\lceil \alpha \rceil}\left[1 + 2X\bar{b}\right].
\end{align*}
%Using previous arguments, the estimation of $I_6$ is given below.
%\begin{align*}
%I_6 \le \epsilon s_0\bar{M}_{\lceil \alpha \rceil} + 2s_0\int_{z_2}^\infty y^\alpha b(x,y)f_n(y,s)\, dy & \le \epsilon\bar{M}_{\lceil \alpha \rceil} + 2s_0\bar{b} \int_{z_2}^\infty y^\alpha f_n(y,s)\, dy \notag \\
%& \le \left(\bar{M}_{\lceil \alpha \rceil} + 2\bar{b}s_0\right)\epsilon. %\label{1.41}
%\end{align*}
The integrand in $I_7$ is continuous in $\bar{\mathcal{S}}$. Therefore, using $(A3)$ and the above arguments we can easily obtain a constant independent of $n$, such that $I_7<\mathcal{B}_2\epsilon.$
Using $\sup_{|x'-x|<\delta} \left|S_n(x')-S_n(x)\right|<\epsilon$ and $(A2)$, $I_8$ and $I_9$ can be estimated as $$I_8 < L\epsilon\quad \mbox{and} \quad I_9 \le s_0X^\alpha \zeta_n(s),$$ respectively.
Now the first integral $I_1$ is written as,
\begin{align*}
I_1 \le \frac{k}{2}(1+X)^\lambda \int_x^{x'} \frac{f_n(x'-y,s)f_n(y,s)}{(x'-y)^\sigma y^\sigma}\, dy & \le \frac{k}{2}(1+X)^\lambda X^\sigma L\int_x^{x'} \frac{f_n(x'-y,s)}{(x'-y)^\sigma}\, dy.
\end{align*}
Substituting $x'-y=z$, we get
\begin{align*}
I_1 \le \frac{k}{2}(1+X)^\lambda X^\sigma L\int_0^{x'-x} \frac{f_n(z,z)}{z^\sigma}\, dz \le \frac{k}{2}(1+X)^\lambda X^{\sigma} L\bar{M}_{-2\sigma}\delta(\epsilon).
\end{align*}
Similarly, for the other integrals $I_2$ and $I_3$ can be obtained as small quantities. Therefore, recalling the relation (\ref{2.22}) and using all the above estimates, we obtain
\begin{align*}
\zeta_n(t) & \le \mathcal{B}_3\epsilon + \mathcal{B}_4\int_0^t \zeta_n(s)\, ds, %\label{2.22}
\end{align*}
where $\mathcal{B}_3$ and $\mathcal{B}_4 (:=\mathcal{B}_1 + s_0X^\alpha)$ are constants independent of $n$ and $\epsilon$. Therefore, applying Gronwall's inequality, we get
\begin{align}
\zeta_n(t)\le \mathcal{B}_3\exp(\mathcal{B}_4 T)\epsilon. \label{2.22}
\end{align}
Hence, the equicontinuity with respect to $x$ on $\bar{\mathcal{S}}$ is obtained.

\subsection{Equicontinuity w.r.t.\ time$-$variable $t$:}
\label{s2.3}
Let $0\le t\le t'\le T$. Therefore, the equation (\ref{2.4}) can be rearranged as follows:
\begin{align}
|f_n(x,t')-f_n(x,t)| \le & \int_t^{t'} \left[\frac{1}{2} \int_0^x K_n(x-y,y)f_n(x-y,s)f_n(y,s)\, dy \right. \notag \\
 & \quad + f_n(x,s)\int_0^\infty K_n(x,y)f_n(y,s)\, dy \notag \\
&\quad + \left. \int_x^\infty b(x,y)S_n(y)f_n(y,s)\, dy + f_n(x,s)S_n(x) \right]\, ds. \label{2.23}
\end{align}
We aim to show that $$\left|f_n(x,t') - f_n(x,t) \right|< \epsilon, \quad \mbox{whenever}\quad \left|t' - t\right| < \delta(\epsilon).$$

We now proceed to estimate the first integral of relation (\ref{2.23}). Without any loss of generality, we consider $\frac{2}{X}\le x\le X$. Therefore, we can write
\begin{align*}
\int_0^x K_n(x-y,y)f_n(x-y,s)f_n(y,s)\, dy \le & k\left( 1 + X \right)^\lambda \int_0^x \frac{1}{(x-y)^\sigma y^\sigma} f_n(x-y,t)f_n(y,t)\, dy \\
= & k\left( 1 + X \right)^\lambda \underbrace{\int_0^{1/X} \frac{f_n(x-y,t)}{(x-y)^\sigma}\frac{f_n(y,t)}{ y^\sigma} \, dy}_{J_4} \\
& + k\left( 1 + X \right)^\lambda \underbrace{\int_{1/X}^x \frac{f_n(x-y,t)}{(x-y)^\sigma}\frac{f_n(y,t)}{ y^\sigma}\, dy}_{J_5}.
\end{align*}
Substituting, $x-y=z$ in $J_4$ and using similar estimate like $I_1$, we get
\begin{align*}
J_4 \le \int_{x-1/X}^x \frac{f_n(z,t)}{z^\sigma}\frac{f_n(x-z,t)}{ (x-z)^\sigma} \, dz \le LX^\sigma\int_{1/X}^x \frac{f_n(x-z,t)}{ (x-z)^\sigma} \, dz < LX^{\sigma}\bar{M}_{-2\sigma}\delta(\epsilon).
\end{align*}
Similarly, $J_5$ can be estimated and hence we obtain
\begin{align*}
\frac{1}{2}\int_0^x K_n(x-y,y)f_n(x-y,s)f_n(y,s)\, dy < kL\left( 1 + X \right)^\lambda X^{\sigma}\bar{M}_{-2\sigma}\delta(\epsilon).
\end{align*}
In this regard, let us denote $\{(x,t):\frac{2}{X}\le x\le X, 0\le t\le T\}=:\tilde{\mathcal{S}}\subset \bar{\mathcal{S}}$.\vspace{0.2cm}\\
For the second integral of the inequality (\ref{2.23}), we recall the estimation of $I_4$ to get a constant $\mathcal{B}_5$ (independent of $n$), such that
$$f_n(x,s)\int_0^\infty K_n(x,y)f_n(y,s)\, dy \le L\mathcal{B}_5.$$
For $(x,t)\in \bar{\mathcal{S}}$, the third integral of (\ref{2.23}) can be estimated as
\begin{align*}
\int_x^\infty S_n(y)b(x,y)f_n(y,s)\, dy = & \frac{1}{x}\int_x^\infty xb(x,y)S_n(y)f_n(y,s)\, dy \\
\le & \bar{b}s_0X\int_0^\infty y^\alpha f_n(y,s)\, dy \le \bar{b}s_0X\bar{M}_{\lceil \alpha \rceil}.
\end{align*}
Finally, the last term of the inequality (\ref{2.23}) can be written as $S_n(x)f_n(x,s)\le Ls_0X^\alpha.$\vspace{0.2cm}\\
Combining the estimates of all the four terms in the relation (\ref{2.23}), we get
\begin{align}
|f_n(x,t')-f_n(x,t)| & < \int_t^{t'} \left[ kL\left( 1 + X \right)^\lambda X^{\sigma}\bar{M}_{-2\sigma}\delta(\epsilon) + L\mathcal{B}_5 + \bar{b}s_0X\bar{M}_{\lceil \alpha \rceil} + Ls_0X^\alpha \right]\, ds\notag \\
& < \mathcal{B}_6\delta(\epsilon), \label{2.24}
\end{align}
where $\mathcal{B}_6 := kL\left( 1 + X \right)^\lambda X^{\sigma}\bar{M}_{-2\sigma}\delta(\epsilon) + L\mathcal{B}_5 + \bar{b}s_0\bar{M}_{\lceil \alpha \rceil} + Ls_0X^\alpha$, is a constant independent of $n$ and $\epsilon$. Therefore, $\{f_n\}$ is equicontinuous with respect to the variable $t$ on $\tilde{\mathcal{S}}$.\vspace{0.2cm}\\
Thus, from the results (\ref{2.22}) and (\ref{2.24}), we can conclude that there exists constants $\mathcal{B}_3$, $\mathcal{B}_4$ and $\mathcal{B}_6$ being independent of $n$ and $\epsilon$, such that
\begin{align}
\sup_{|x'-x|<\delta,~|t'-t|<\delta}\left|f_n(x',t')-f_n(x,t)\right|\le \left[\mathcal{B}_3\exp(\mathcal{B}_4 T) + \mathcal{B}_6\right]\epsilon, \label{2.25}
\end{align}
whenever $\frac{2}{X}\le x\le x' \le X$, $0\le t\le t' \le T$.
\end{proof}
Furthermore, combining the relations (\ref{2.15}), (\ref{2.25}) along with the Arzel\`{a}-Ascoli theorem [\cite{ash1972measure,edwards1994functional}], we get the existence of a subsequence $\left\{ f_{n_s}\right\}_{s=1}^\infty$, which is relatively compact in the uniform-convergence topology of continuous functions on each rectangle $\tilde{\mathcal{S}}$.

\subsection{Proof of Theorem \ref{th1}:}
\label{s2.4}
By means of diagonal method we can select a subsequence $\left\{ f_{p}\right\}_{p=1}^\infty$ of the sequence $\left\lbrace f_n \right\rbrace_{n=1}^\infty$ that converges uniformly on each compact subset of $\mathcal{S}$ to a continuous non-negative function $f$. Let us consider an integral $$\int_{\bar{z}_1}^{\bar{z}_2} \left(x^{j}+ \frac{1}{x^{2\sigma}}\right) f(x,t)\, dx, \quad 0\le j\le r,~\sigma\ge 0.$$
Now, the subsequence $\{f_p\}$ ensures that for all $\epsilon>0$ there exists $p\ge 1$ such that, the following relation holds good.
\begin{align}
\int_{\bar{z}_1}^{\bar{z}_2} \left(x^{j} + \frac{1}{x^{2\sigma}}\right)f(x,t)\, dx\le \int_{\bar{z}_1}^{\bar{z}_2} \left(x^{j} + \frac{1}{x^{2\sigma}}\right) f_{p}(x,t)\, dx + \epsilon. \label{2.26}
\end{align}
Furthermore, $\epsilon$, $\bar{z}_1$ and $\bar{z}_2$ are arbitrary in (\ref{2.26}). Therefore, for all $0\le j\le r$ and $2\sigma \ge 0$, we have
\begin{align}
\int_0^\infty \left(x^{j} + \frac{1}{x^{2\sigma}}\right) f(x,t)\, dx \le \bar{M}_{j} + \bar{M}_{-2\sigma}. \label{2.27}
\end{align}
We aim to show that $f(x,t)$ is a solution to the initial-value problem (\ref{1.1}) and (\ref{1.2}). In this regard, let us replace $K_n$, $S_n$, $f_n$ in (\ref{2.4}) with $K_p-K+K$, $S_p-S+S$, $f_p-f+f$ respectively, and rearrange the terms to obtain
\begin{align}
(f_p-f)(x,t)+f(x,t)= f_0(x) & + \int_0^t \left[ \int_0^x \frac{1}{2}(K_p-K)(x-y,y)f_p(x-y,s)f_p(y,s)\, dy \right. \notag \\
& + \int_0^x \frac{1}{2}K(x-y,y)(f_p-f)(x-y,s)f_p(y,s)\, dy \notag \\
& + \int_0^x \frac{1}{2}K(x-y,y)(f_p-f)(y,s)f(x-y,s)\, dy \notag \\
& + \int_0^x \frac{1}{2}K(x-y,y)f(x-y,s)f(y,s)\, dy \notag\\
& - f_p(x,s)\int_0^\infty (K_p-K)(x,y)f_p(y,s)\, dy \notag\\
& - (f_p-f)(x,s)\int_0^\infty K(x,y)f_p(y,s)\, dy \notag \\
& - f(x,s) \int_0^\infty K(x,y)(f_p-f)(y,s)\, dy \notag\\
& - f(x,s) \int_0^\infty K(x,y)f(y,s)\, dy \notag \\
& + \int_x^\infty b(x,y)(S_p-S)(y)f_p(y,s)\, dy  \notag \\
& + \int_x^\infty b(x,y)S(y)(f_p-f)(y,s)\, dy \notag \\
& + \left.\int_x^\infty b(x,y)S(y)f(y,s)\, dy - (S_p-S)(x)f_p(x,s)\right] \, ds\notag \\
& - \int_0^t \left[ S(x)(f_p-f)(x,s) - S(x)f(x,s)\right]\,ds. \label{2.28}
\end{align}
Using previous arguments, we can easily get
\begin{align*}
\int_0^\infty \left|K_p-K\right|(x,y)f_p(y,s)\, dy  \le \mathcal{B}_7\epsilon, & \quad \int_0^\infty K(x,y)\left|f_p-f\right|(y,s)\, dy \le \mathcal{B}_8\epsilon,\\
\int_x^\infty b(x,y)\left|S_p-S\right|(y)f_p(y,s)\, dy \le \mathcal{B}_9\epsilon, & \quad \int_x^\infty b(x,y)S(y)\left|f_p-f\right|(y,s)\, dy \le \mathcal{B}_{10}\epsilon
\end{align*}
where $\mathcal{B}_7$, $\mathcal{B}_8$, $\mathcal{B}_9$ and $\mathcal{B}_{10}$ are constants. Hence, in (\ref{2.28}) when $p\rightarrow \infty$, all the terms involving integrals over the infinite range tend to zero due to the estimates of their \emph{tails}. Furthermore, it has already been proved earlier that the difference terms of equation (\ref{2.28}) involving the integrals over finite ranges are convergent. Hence, using the definition of convergence of a sequence, it can easily be proved that all the integrals involving $(K_p-K)$, $(S_p-S)$ and $(f_p-f)$ tend to zero as $p\rightarrow \infty$. Thus, in the limiting case (\ref{2.28}) reduces to
\begin{align}
f(x,t) = f_0(x) + \int_0^t & \left[ \frac{1}{2}\int_0^x K(x-y,y)f(x-y,s)f(y,s)\, dy - f(x,s)\int_0^\infty K(x,y)f(y,s)\, dy \right. \notag\\
& \left. + \int_x^\infty b(x,y)S(x)f(y,s)\, dy - S(x)f(x,s) \right]\, ds. \label{2.29}
\end{align}
The estimates of the tails and the continuity of $f(x,t)$ together implies that the r.h.s.\ of (\ref{1.1}), evaluated at $f$, is a continuous function on $\mathcal{S}$. Moreover, relation (\ref{2.26}) ensures that  $f(x,t)\in \Upsilon_{r,\sigma}^+(T)$ and differentiation of (\ref{2.29}) with respect to $t$ establishes that $f(x,t)$ is a continuous and differentiable solution of the IVP (\ref{1.1}), (\ref{1.2}).
\end{proof}

\section{Uniqueness theory}
\label{s3}

In this section, we prove the uniqueness of the solutions in the space $\tilde{\Upsilon}_{r,\sigma}^+(T)$, where $r\ge 2$, $\sigma\ge 1$ and $\tilde{\Upsilon}_{r,\sigma}^+(T)$ denotes the positive cone of the Banach space defined by
$$\tilde{\Upsilon}_{r,\sigma} (T) = \left\lbrace f(x,t)\in \mathcal{C}\left(]0,\infty[,[0,T] \right) : \left\|f\right\|_{\tilde{\Upsilon}} < \infty \right\rbrace$$ endowed with the norm $$\left\|f\right\|_{\tilde{\Upsilon}}:=\sup_{0\le t\le T} \int_0^\infty \left( x^{r} + \frac{1}{x^{2\sigma}}\right) \left|f(x,t)\right|\, dx.$$ We now state and prove the following theorem:

\begin{Theorem}
\label{th2}
Let the functions $K(x,y)$, $b(x,y)$ and $S(x)$ be non-negative and continuous for all $x,y\in]0,\infty[$ and satisfy the condition $(A1)$, $(A2)$ and $(A3)$. If the initial data $f_0(x)$ belongs to $\tilde{\Upsilon}_{r,\sigma}^+(0)$, then the IVP (\ref{1.1}), (\ref{1.2}) has a unique solution in $\tilde{\Upsilon}_{r,\sigma}^+(T)$.
\end{Theorem}

\begin{proof}
Let for $t\ne 0$, $v_1(x,t)$ and $v_2(x,t)$ be two distinct solutions of (\ref{1.1}), (\ref{1.2}) along with $v_1(x,0) = v_2(x,0)$. Let $V(x,t) := v_1(x,t)-v_2(x,t)$ and we construct an auxiliary function $$\mathcal{V}(t) := \int_0^\infty \left( x^2 + \frac{1}{x^{\sigma}} \right)\left|V(x,t)\right| \, dx.$$ Following the construction of $V(x,t)$, we get it to be absolutely continuous over the interval $]0,\infty[\times[0,T]$ and hence the function $V(x,t)$ satisfies equation (\ref{1.1}). Therefore, the derivative of the solutions is obtained as
\begin{align}
\frac{\partial V(x,t)}{\partial t} = & \frac{1}{2}\int_0^x K(x-y,y)\left[ v_1(x-y,t)v_1(y,t)-v_2(x-y,t)v_2(y,t)\right]\, dy \notag \\
& - \int_0^\infty K(x,y)\left[ v_1(x,t)v_1(y,t)-v_2(x,t)v_2(y,t) \right]\, dy \notag \\
& + \int_x^\infty b(x,y)S(y)V(y,t)\, dy - S(x)V(x,t). \label{3.1}
\end{align}
For all $t\in\R$, we define
$$\sgn(t)=\left \{ \begin{array}{ll}
    ~1, &\mbox{when $ t>0, $}\notag \\
    ~0, &\mbox{when $ t=0, $}\notag \\
    -1, &\mbox{when $t<0, $}
     \end{array}
\right. \quad \mbox{and} \quad \frac{d|w(t)|}{dt} = \sgn(w(t))\frac{dw(t)}{dt}.$$
Multiplying both sides of (\ref{3.1}) by the weight $\left( x^2 + \frac{1}{x^{\sigma}}\right)$ then integrating with respect to $x$ over the range $]0,\infty[$ and using the fact that $\mathcal{V}(0)=0$, we get
\begin{align}
\mathcal{V}(t) = & \int_0^t \int_0^\infty \left( x^2 + \frac{1}{x^{\sigma}}\right)\sgn\left(V(x,s)\right) \notag \\
& \times \left[\underbrace{\frac{1}{2}\int_0^x K(x-y,y)\left[ v_1(x-y,s)v_1(y,s)-v_2(x-y,s)v_2(y,s)\right]\, dy}_{\mathcal{V}_1} \right. \notag \\
& \quad - \int_0^\infty K(x,y)\left[ v_1(x,s)v_1(y,s)-v_2(x,s)v_2(y,s) \right]\, dy \notag \\
& \quad + \left. \underbrace{\int_x^\infty b(x,y)S(y)V(y,s)\, dy}_{\mathcal{V}_2} - \underbrace{S(x)V(x,s)}_{\mathcal{V}_3}\right]\, dx\, ds. \label{3.2}
\end{align}
Let $$\mathcal{V}_4 := K(x,y)\left[ v_1(x,s)v_1(y,s)-v_2(x,s)v_2(y,s) \right].$$ Considering the integral $\mathcal{V}_1$, changing the order of integration and substituting $x-y=x'$, $y=y'$, we obtain
\begin{align*}
&\frac{1}{2}\int_0^\infty \left( x^2 + \frac{1}{x^{\sigma}}\right)\sgn\left(V(x,s)\right)\mathcal{V}_1\, dx \\
= &\frac{1}{2} \int_0^\infty \int_0^\infty \left( (x+y)^2 + \frac{1}{(x+y)^{\sigma}}\right)\sgn\left(V(x+y,s)\right) \mathcal{V}_4\, dx\, dy. %\label{3.6}
\end{align*}
Using the above relation and the symmetry of $K(x,y)$, the term $\mathcal{V}(t)$ is written as
\begin{align}
\mathcal{V}(t) =  & \frac{1}{2}\int_0^t \int_0^\infty \int_0^\infty \left[ \left((x+y)^2 + \frac{1}{(x+y)^{\sigma}}\right) \sgn\left(V(x+y,s)\right) \right. \notag \\
& \quad \quad \quad \quad \left. - \left(x^2 + \frac{1}{x^{\sigma}}\right)\sgn\left(V(x,s)\right) - \left( y^2 + \frac{1}{y^{\sigma}}\right)\sgn\left(V(y,s)\right)\right]\mathcal{V}_4 \, dx \, dy\, ds\notag \\
& + \int_0^t \int_0^\infty \left( x^2 + \frac{1}{x^{\sigma}}\right)\sgn\left(V(x,s)\right)\left[\mathcal{V}_2-\mathcal{V}_3 \right]\, dx\, ds. %\label{3.13}
\end{align}
Now, for $x,y\ge 0$ and $s\in [0,T]$, we define the following weight function
\begin{align*}
w(x,y) := & \left[\left((x+y)^2 + \frac{1}{(x+y)^{\sigma}}\right)\sgn\left(V(x+y,s)\right) - \left(x^2 + \frac{1}{x^{\sigma}}\right)\sgn\left(V(x,s)\right)\right. \\
& \quad  \left. - \left(y^2 + \frac{1}{y^{\sigma}}\right)\sgn\left(V(y,s)\right)\right].
\end{align*}
Further, we have the following relation
\begin{align*}
v_1(x,s)v_1(y,s) - v_2(x,s)v_2(y,s) = & v_1(x,s)(v_1-v_2)(y,s) + v_2(y,s)(v_1-v_2)(x,s)\\
= & v_1(x,s)V(y,s) + v_2(y,s)V(x,s).
\end{align*}
Therefore, $\mathcal{V}(t)$ can be represented as
\begin{align*}
\mathcal{V}(t) = & \frac{1}{2}\int_0^t\int_0^\infty\int_0^\infty w(x,y)K(x,y)[v_1(x,s)V(y,s) + v_2(y,s)V(x,s)]\, dy\, dx\,ds \notag \\
& + \int_0^t\int_0^\infty \left( x^2 + \frac{1}{x^{\sigma}}\right)\sgn\left(V(x,s)\right)\left[ \mathcal{V}_2-\mathcal{V}_3\right] \, dx\,ds \notag \\
= & \underbrace{\frac{1}{2}\int_0^t\int_0^\infty\int_0^\infty w(x,y)K(x,y)v_1(x,s)V(y,s)\, dy\, dx\,ds}_{I_a} \\%\label{3.14}\\
& + \underbrace{\frac{1}{2}\int_0^t\int_0^\infty\int_0^\infty w(x,y)K(x,y)v_2(y,s)V(x,s)\, dy\, dx\,ds}_{I_b} \\%\label{3.15}\\
& + \underbrace{\int_0^t\int_0^\infty \left( x^2 + \frac{1}{x^{\sigma}}\right)\sgn\left(V(x,s)\right)\left[ \mathcal{V}_2-\mathcal{V}_3 \right]\, dx\,ds}_{I_c}. %\label{3.16}
\end{align*}
Now, for all $x,y>0$, we have
$$\frac{1}{(x+y)^p}\le \frac{1}{x^p} + \frac{1}{y^p},\mbox{ if } p>0 \quad \mbox{and} \quad (x+y)^p \le x^p +y^p, \mbox{ if } 0<p<1,$$ and taking in account that for all $t_1,t_2\in\R$ that $\sgn(t_1)\sgn(t_2)=\sgn(t_1t_2)$ and $|t_1|=t_1\sgn(t_1)$, we obtain
\begin{align}
w(x,y)V(y,t) \le & \left[ \left((x+y)^2 + \frac{1}{(x+y)^{\sigma}}\right) - \left(x^2 + \frac{1}{x^{\sigma}}\right) - \left( y^2 + \frac{1}{y^{\sigma}}\right)\right]|V(y,t)| \notag \\
\le & 2xy|V(y,s)|. \label{3.18}
\end{align}
So, the integral $I_a$ gives
\begin{align}
I_a\le & k\int_0^t\underbrace{\int_0^\infty x^{1-\sigma}(1+x)^\lambda v_1(x,s)\, dx}_{I_{a_1}}\underbrace{\int_0^\infty y^{1-\sigma}(1+y)^\lambda|V(y,s)|\, dy}_{I_{a_2}}\, ds.
\end{align}
Now, estimating $I_{a_1}$
\begin{align*}
I_{a_1} = & \int_0^1 x^{1-\sigma}(1+x)^\lambda v_1(x,s)\, dx + \int_1^\infty x^{1-\sigma}(1+x)^\lambda v_1(x,s)\, dx\\
\le & 2\left[ \int_0^1 \left(x^{1-\sigma} + x^{\lambda-\sigma+1}\right) v_1(x,s)\, dx + \int_1^\infty \left(x^{1-\sigma} + x^{\lambda-\sigma+1}\right) v_1(x,s)\, dx\right],
\end{align*}
\begin{enumerate}
\item when $x\le 1$ then both of $x^{1-\sigma}$, $x^{\lambda-\sigma+1} \le 1$,
\item when $x\ge 1$ then both of $x^{1-\sigma}$, $x^{\lambda-\sigma+1} \le x^2$ as $\lambda-\sigma+1\in[1,2]$.
\end{enumerate}
Therefore,
\begin{align*}
I_{a_1} \le & 4\left[ \int_0^1 v_1(x,s)\, dx + \int_1^\infty x^2 v_1(x,s)\, dx \right] \\
\le & 4\left[ \int_0^1 \frac{1}{x^{2\sigma}}v_1(x,s)\, dx + \int_1^\infty x^{r} v_1(x,s)\, dx \right] ,\quad [x\le 1\Rightarrow \frac{1}{x}\ge 1.]\\
\le & 4\left[ \int_0^1 \left( x^{r} + \frac{1}{x^{2\sigma}}\right) v_1(x,s)\, dx + \int_1^\infty \left(x^{r} + \frac{1}{x^{2\sigma}}\right) v_1(x,s)\, dx \right] \\
\le & 4\int_0^\infty \left(x^{r}+\frac{1}{x^{2\sigma}}\right)v_1(x,s)\, dx = 4\|v_1\|_{\tilde{\Upsilon}}.
\end{align*}
Similarly, it can be shown that $I_{a_2} \le 4\mathcal{V}(s).$
Therefore, $$I_a \le ~\Gamma_1\int_0^t \mathcal{V}(s)\, ds,\quad \mbox{where} \quad \Gamma_1 := 16k\|v_1\|_{\tilde{\Upsilon}}.$$
In a similar way we can show that $I_b \le \Gamma_2\int_0^t \mathcal{V}(s)\, ds.$ \vspace{0.2cm}\\
The third integral $I_c$ is estimated as follows,
\begin{align*}
I_c = & \underbrace{\int_0^t\int_0^\infty x^2\left[ \int_x^\infty b(x,y)S(y)|V(y,s)|\, dy - S(x)|V(x,s)| \right] \, dx\,ds}_{I_{c_1}} \\
& + \underbrace{\int_0^t\int_0^\infty \frac{1}{x^{\sigma}}\left[ \int_x^\infty b(x,y)S(y)|V(y,s)|\, dy - S(x)|V(x,s)| \right] \, dx\,ds}_{I_{c_2}}.
\end{align*}
Changing the order of integrations of the first terms of both $I_{c_1}$ and $I_{c_2}$, we get
\begin{align*}
I_{c_1} = \int_0^t\left[ \int_0^\infty \int_0^y x^2b(x,y)S(y)|V(y,s)|\, dx\, dy - \int_0^\infty y^2S(y)|V(y,s)|\, dy \right]\,ds,
\end{align*}
\begin{align*}
\mbox{and,} \quad I_{c_2} = \int_0^t\int_0^\infty \left[ \int_0^y\frac{1}{x^{\sigma}} b(x,y)S(y)|V(y,s)|\, dx - S(y)|V(y,s)| \right] \, dy\,ds.
\end{align*}
Considering $I_{c_1}$, we have $x\le y$ and use the relation (\ref{1.3}), to get
\begin{align*}
I_{c_1} \le \int_0^t \left[\int_0^\infty y^2S(y)|V(y,s)|\, dy - \int_0^\infty y^2S(y)|V(y,s)|\, dy\right]\, ds = 0.
\end{align*}
For the relation $I_{c_2}$, using the positivity of $S(y)$ and fact $x\le y$ implies $\frac{y^\sigma}{x^\sigma}\ge 1$ along with the assumption $(A3)$, we get
\begin{align*}
I_{c_2} \le \int_0^t \int_0^\infty \int_0^y \frac{1}{x^{2\sigma}}b(x,y)y^{\sigma}S(y)|V(y,s)|\, dx \le & \bar{N}\int_0^t\int_0^\infty \frac{1}{y^{\sigma}}S(y)|V(y,s)|\, dy\,ds \\
\le & \bar{N}s_0\int_0^t\int_0^\infty y^{\alpha-\sigma}|V(y,s)|\, dy\,ds \\
\le & ~\bar{N}s_0\int_0^t \mathcal{V}(s)\, ds.
\end{align*}
Thus $\mathcal{V}(t)$ is estimated as
\begin{align*}
\mathcal{V}(t) \le ~\Gamma \int_0^t \mathcal{V}(s) \, ds, \quad \mbox{ where } \Gamma:=\left[ \Gamma_1 +\Gamma_2 +\bar{N}s_0 \right].
\end{align*}
Applying Gronwall's Lemma, we get
\begin{align*}
\mathcal{V}(t) = 0,\quad \mbox{which implies,}\quad v_1(x,t) = v_2(x,t) \quad \mbox{for all}\quad 0\le t\le T.
\end{align*}
Hence, the uniqueness of the solution is obtained.
\end{proof}
%\begin{Remark}
%Thus in this section we have established the uniqueness of the whole class of strong solutions of the coagulation fragmentation equation \eqref{1.1} \& \eqref{1.2} having singularity in its coagulation kernel. The domain of uniqueness being a restricted part of the domain of existence. In the previous works of \cite{camejo2014singular,camejo2014regular}, the authors used to prove uniqueness property for a subclass of the solutions. So this can be treated as a interesting achievement of this part of our study.
%\end{Remark}

\section{Conclusion}
\label{s5}
In this work, a thorough mathematical investigation on the existence and uniqueness of strong solutions to the continuous CF-equation (\ref{1.1}), (\ref{1.2}) is performed. Here, the coagulation kernel includes a huge class of singular functions. Therefore, we suitably construct a sequence of truncated problems, in order to take care the singular coagulation kernel. Using some standard results from \cite{dubovskii1994mathematical}, the existence of solutions for these truncated problems is obtained. Next, we establish that the solutions are uniformly bounded and equicontinuous both in $x$ and $t$, over a compact rectangular strip. Now, using the Arzel\`a-Ascoli theorem, a strongly convergent subsequence of the solutions is obtained, whose limiting function is proved to be the solution of the IVP (\ref{1.1}), (\ref{1.2}). The study is completed by establishing the uniqueness of the solutions under some additional restrictions. Beside mathematical significance, the present study includes several physically important singular coagulation kernels like, Brownian diffusion kernel, granulation kernel, EKE kernel, activated sludge flocculation kernel etc.\ [Please see Appendix-\ref{Apdx1}, for further details]. However, the non-random coalescence kernel of \cite{sastry1975similarity} is still beyond the scope of our present study.

\section{Acknowledgement}
Author JS thank Deutscher Akademischer Austausch Dienst (DAAD) and Ministry of Human Resource Development (MHRD), Government of India for their funding support. Author JK thank Alexander von Humboldt (AvH) Foundation for his funding support. Both the authors gratefully acknowledges Otto-von-Guericke University, Magdeburg, Germany for their hospitality during this work.

\titleformat{\section}{\large\bfseries}{\appendixname~\thesection .}{0.5em}{}
\begin{appendices}
\section{Singular kernels}
\label{Apdx1}
\begin{enumerate}
\item \emph{Brownian diffusion kernel}: $K(x,y) \le k\left( x^{1/3} + y^{1/3} \right)\left( x^{-1/3} + y^{-1/3} \right)$ [\cite{smoluchowski1917}].
\item \emph{Granulation kernel}: $K(x,y) \le k(x+y)^a(xy)^{-b}$, $a,b>0$ [\cite{kapur1972kinetics}].
\item \emph{Equi-partition of kinetic energy (EKE) kernel}: $K(x,y) \le k\left( x^{1/3} + y^{1/3} \right)^2\sqrt{\frac{1}{x}+\frac{1}{y}}$ [\cite{hounslow1998population}].
\item \emph{Activated sludge flocculation kernel}: $k\left( x^{1/3} + y^{1/3} \right)^q\left[8y_c + \left(x^{1/3} + y^{1/3}\right)^3\right]^{-1},\quad 0\le q\le3$ and $y_c$ is a critical volume [\cite{Ding}].
\item \emph{Non-random coalescence kernel}: $K(x,y) \le k(x^{2/3}+y^{2/3})\left(1/x +1/y\right)$ [\cite{sastry1975similarity}].
\end{enumerate}
\end{appendices}
%\end{spacing}
%\section{References}
%\bibliographystyle{agsm}\small
\bibliographystyle{harvard}\small
\bibliography{existence}
\end{document}